\documentclass[12pt]{amsart}

\usepackage{epsfig}
\usepackage{amssymb}
\usepackage{latexsym}
\usepackage{comment}

\newtheorem{theorem}{Theorem}[section]

\newtheorem{proposition}[theorem]{Proposition}

\numberwithin{figure}{section} \numberwithin{equation}{section}

\title{Global Dynamics of the smallest Chemical Reaction System with
 Hopf Bifurcation}

\author{
{  Hal L. Smith*\\
School of Mathematical and Statistical Sciences\\
 Arizona State University\\
  Tempe, AZ, USA 85287 }\\}

\thanks{*Supported by NSF Grant DMS-0918440}

\begin{document}

\maketitle

\begin{abstract}
The global behavior of solutions is described for the smallest chemical reaction system that exhibits a Hopf bifurcation, discovered in \cite{WH1}. This three-dimensional system is a competitive system and a monotone cyclic feedback system. The Poincar\'e-Bendixson theory extends to such systems \cite{MS,H0,HS,S} and a Bendixson criterion exists to rule out periodic orbits \cite{LM}.
\end{abstract}

\pagestyle{myheadings}\markboth{\sc  H.L. Smith} {\sc
Global Dynamics of a Reaction System}


\section{Introduction and Main Results}

 In \cite{WH1}, Wilhelm and Heinrich  discover the ``smallest chemical reaction system''  which may exhibit a Hopf bifurcation. Like the famous Lorenz system, it is a three dimensional system with only a single quadratic nonlinearity. In \cite{WH2}, they study the Hopf bifurcation  in detail using center manifold and normal form techniques which establish the existence of
 a stable periodic orbit very near the bifurcation point. However, as they point out, their analysis does not preclude that chaotic behavior may not occur in this simple system. Sprott \cite{Sp} has cataloged many other three dimensional systems with quadratic nonlinearities which can have chaotic dynamics. However, we will show that this particular chemical reaction system does not support chaos. We  give a fairly complete
 analysis of the global dynamics of the system of ordinary differential equations studied in \cite{WH2}. It is given by the following equations:

\begin{eqnarray}
  x' &=& kx-k_2xy \nonumber\\
  y' &=& k_5z-k_3y \\
  z' &=& k_4 x-k_5z \nonumber
\end{eqnarray}
where $k_i>0,\ i=2,3,4,5$ and $k=k_1 A-k_4$ where $A$ is the concentration of ``the outer reactant of the autocatalytic reaction''. The $k_i$ are rate constants and $k$ need not be positive; see \cite{WH1}. Variables $x,y,z$ denote concentrations and therefore
are nonnegative.

The scaling $\bar x=x/a,\ \bar y=y/b,\ \bar z=z/c$ and the choice
$k_2 b=1,\ k_5c=k_3b,\ k_4 a=k_5 c$ leads to the system:
\begin{eqnarray}\label{scaled}
  x' &=& kx-xy \nonumber\\
  y' &=& k_3(z-y) \\
  z' &=& k_5(x-z)\nonumber
\end{eqnarray}
where we have dropped the bar over variables for simplicity. Our focus is on the dynamics exhibited by \eqref{scaled} in
the nonnegative octant $\mathbb{R}^3_+$.

System \eqref{scaled} is a competitive system because the signed incidence matrix of its Jacobian has a single loop with
two positive feedbacks and one negative one:
$$
x \xrightarrow{+} z \xrightarrow{+} y \xrightarrow{-} x
$$
See \cite{S,HS}; the change of variables $z\to-z$ gives the system the canonical
form with all off-diagonal entries of the Jacobian being non-positive. Since it is three dimensional, the celebrated Poincar\'e-Bendixon
Theorem extends to its solutions by a result of M.W. Hirsch \cite{H0,S,HS}: A compact omega or alpha limit set that contains no equilibrium is a periodic orbit.
In fact, system \eqref{scaled} is also a monotone cyclic feedback system \cite{MS}, for which the Poincar\'e-Bendixon Theorem also holds (regardless
of the dimension of the system).

We note that $0$ is the unique equilibrium of \eqref{scaled} if $k\le 0$; it is asymptotically stable if $k\le 0$ and unstable when $k>0$.
In the latter case, the stable manifold of $0$ in $\mathbb{R}^3_+$ is $S=\{(x,y,x)\in \mathbb{R}^3_+:x=0\}$.
If $k>0$, there is an additional equilibrium $E=k(1,1,1)$ which is asymptotically stable if $k<k_3+k_5$ and unstable if $k>k_3+k_5$ with a pair
of complex conjugate eigenvalues with positive real part and one negative eigenvalue.

Part (a) of the following result was proved in \cite{WH2} using a linear Lyapunov function.

\begin{theorem}\label{solnbehavior}
The following hold for \eqref{scaled}:
\begin{enumerate}
  \item [(a)] If $k\le 0$, every solution converges to $0$.
  \item [(b)] If $0<k<k_3+k_5$, then every solution except those on $S$ converge to $E$.
    \item [(c)] If $k>k_3+k_5$, then every solution not starting on $S$ or on the one-dimensional stable manifold
    of $E$, converges to a nontrivial periodic orbit. At least one periodic orbit is orbitally asymptotically stable and the number of
    periodic orbits is finite.
\end{enumerate}
Moreover, if $k>0$ then any solution not starting in $S$ satisfies
\begin{equation}\label{timeave}
\lim_{T\to\infty}\frac{1}{T}\int_0^T u(s)ds=k,\quad u=x,y,z.
\end{equation}
\end{theorem}

The main open problem concerning the dynamics of \eqref{scaled} is the number of periodic orbits in case (c). Numerical
simulations suggest there is only one.

The geometry of the stable manifold of $E$ is described in the following result. It is tangent at $E$ to the eigenvector
with sign pattern $(+,+,-)$ and has certain monotonicity properties due to the monotonicity of the time-reversed system.

\begin{proposition}\label{stablemanE}
Let $k>k_3+k_5$. Then the stable manifold of $E$ in $\mathbb{R}^3_+$ consists of $E$ and two monotone solutions
$p_u(t)=(x_u(t),y_u(t),z_u(t)), \ t\in [0,\infty)$ and $p_l(t)=(x_l(t),y_l(t),z_l(t)), \ t\in [0,\infty)$.

$p_u$ satisfies $p_u(0)=(x_u(0),0,z_u(0))$ with $0<x_u(0)<k<z_u(0)$, $p_u(\infty)=E$ and $x_u,y_u$ are strictly increasing while $z_u$ is strictly decreasing. Therefore, the graph of $p_u$ belongs to $[0,k]\times[0,k]\times[k,\infty)$.

$p_l$ satisfies $p_l(0)=(x_l(0),y_l(0),0)$ with $x_l(0),y_l(0)>k$, $p_l(\infty)=E$ and $x_l,y_l$ are strictly decreasing while $z_l$ is strictly increasing. Therefore, the graph of $p_l$ belongs to $[k,\infty)\times [k,\infty)\times [0,k]$.

\end{proposition}

As a consequence of Proposition~\ref{stablemanE}, the unstable manifold of $0$ connects $0$ to a nontrivial
periodic solution in case (c) of Theorem~\ref{solnbehavior}; it cannot be a heteroclinic orbit connecting $0$ to $E$.

\begin{theorem}\label{attractor}
System \eqref{scaled} has a non-empty compact invariant set $A$ that attracts all bounded subsets of $\mathbb{R}^3_+$.

If $k\le 0$, then $A=\{0\}$.

If $k>0$. Then $A\subset [0,M]^3$
where $M=\frac{k(k+k_5)(k+k_3)}{k_3k_5}$.
\end{theorem}

\begin{theorem}\label{persistence}
If $k>0$, then \eqref{scaled} is uniformly persistent, i.e., $\exists \eta>0$
such that:
\begin{equation}\label{persist}
\min\{x(0),y(0),z(0)\}>0 \Rightarrow \liminf_{t\to\infty} \min\{x(t),y(t),z(t)\}>\eta.
\end{equation}

The attractor $A$ is the disjoint union $A=\{0\}\cup C\cup A_1$ where:
\begin{enumerate}
  \item [(i)] $\{0\}$ attracts all bounded subsets of $S$.
  \item [(ii)] $A_1$ is a compact invariant subset of the interior of $\mathbb{R}^3_+$ that attracts all compact
  subsets of $\mathbb{R}^3_+$ that do not intersect $S$.
  \item [(iii)] $C$ is the one-dimensional unstable manifold of $\{0\}$ which connects $0$ to $A_1$.
\end{enumerate}

\end{theorem}

\section{Proofs}

\begin{proof} Proof of Theorem~\ref{attractor}.
Let $\phi(t,p)$ denote the solution of \eqref{scaled} which at $t=0$ is at $p=(x(0),y(0),z(0))$.
If $k\le 0$, then $x'\le 0$ and it is trivial to show that $A=\{0\}$.

Assume that $k>0$. According to Theorem 2.33 in \cite{ST}, we must show that \eqref{scaled} is point dissipative and for every bounded set $B$ there exists $t_r\ge 0$
such that $\phi([t_r,\infty)\times B)$ is bounded.

We start by showing that it is point dissipative.
If $x(0)=0$, then $x(t)=0$ for all $t$ and there exists $\eta,\alpha>0$ such that
$y(t)+z(t)\le \eta (y(0)+z(0))e^{-\alpha t}$. Assume that $x(0)>0$. Then direct computation
shows that
\begin{eqnarray}\label{inequal}
  (y/x)' &\ge& k_3(z/x)-(k+k_3)(y/x) \\
  (z/x)' &\ge& k_5-(k+k_5)(z/x)\nonumber
\end{eqnarray}
It follows that $(z/x)_\infty\ge \frac{k_5}{k+k_5}$ and $(y/x)_\infty\ge \frac{k_3k_5}{(k+k_3)(k+k_5)}$
where $f_\infty=\liminf_{t\to\infty}f(t)$. Consequently, for every $\epsilon>0$, there exists $t_0>0$ such that
$$
x'=kx-x^2(y/x)\le kx-m_\epsilon x^2,\quad t>t_0
$$
where $m_\epsilon=\frac{k_3k_5}{(k+k_3)(k+k_5)}-\epsilon$. Hence, $x^\infty=\limsup_{t\to\infty}x(t)\le k/m_\epsilon$
and, since $\epsilon>0$ is arbitrary, $x^\infty\le M$. Using this estimate and the differential equation for $z$ immediately yields
$z^\infty\le M$; using this last estimate and the differential equation for $y$ yields $y^\infty\le M$. Thus, we have shown that all points are attracted to the bounded set $[0,M]^3$ and \eqref{scaled} is point dissipative.

Now we identify a family of positively invariant bounded sets. Let $0<\sigma \le\frac{k_5}{k+k_5}$, $0<\rho\le \frac{k_3\sigma}{k+k_3}$,
and $K\ge k/\rho$. Define
$$
B_{(\sigma,\rho,K)}=\{(x,y,x)\in [0,K]^3:x=0 \ \hbox{or}\ z/x\ge \sigma\ \hbox{and}\ \ y/x\ge \rho\}
$$
We claim that $B_{(\sigma,\rho,K)}$ is positively invariant. Indeed, integrating the second inequality \eqref{inequal}
leads to
\begin{equation}\label{z/x}
(z/x)(t)\ge \frac{k_5}{k+k_5}(1-e^{-(k+k_5)t})+e^{-(k+k_5)t} (z/x)(0)\ge \sigma
\end{equation}
for any solution starting in $B_{(\sigma,\rho,K)}$ with $x(0)\ne 0$. Similarly, using that $(z/x)(t)\ge \sigma$ in the first inequality \eqref{inequal} leads to
$$
(y/x)(t)\ge \frac{k_3\sigma}{k+k_3}(1-e^{-(k+k_3)t})+e^{-(k+k_3)t}(y/x)(0)\ge \rho.
$$
As $x'=kx-x^2(y/x)\le x(k-\rho x)\le 0$ when $x=K$, $z'=k_5(x-z)\le 0$ when $z=K$ and $0\le x\le K$, and
$y'=k_3(z-y)\le 0$ when $y=K$ and $0\le z\le K$, the positive invariance of $B_{(\sigma,\rho,K)}$ follows.

Finally, we show that  $\forall L>0$ there exists $0<\sigma <\frac{k_5}{k+k_5}$, $0<\rho<\frac{k_3\sigma}{k+k_3}$,
and $K\ge k/\rho$ such that $\phi(2, [0,L]^3)\subset B_{(\sigma,\rho,K)}$. Since $B_{(\sigma,\rho,K)}$ is bounded and positively invariant, this proves that  $\phi([2,\infty)\times [0,L]^3)$ is bounded. Hence, bounded sets have uniformly bounded orbits.

Let $(x(0),y(0),z(0))\in [0,L]^3$ and $x(0)>0$.
Integrating the second of \eqref{inequal} gives
$$
\frac{z(t)}{x(t)}\ge \frac{k_5}{k+k_5}(1-e^{-(k+k_5)}),\ t\ge 1.
$$
Let $\sigma$ be the right side of the above inequality and note that $\sigma<\frac{k_5}{k+k_5}$.
Substituting this estimate
into the first of \eqref{inequal}
and integrating leads to
$$
\frac{y(2)}{x(2)}\ge \frac{k_3\sigma}{k+k_3}(1-e^{-(k+k_3)})
$$
Let $\rho$ be the right hand side of the above inequality and note that $\rho<\frac{k_3\sigma}{k+k_3}$.
Now we can choose $K>k/\rho$ so large that $x(2),y(2),z(2)\le K$ holds for every point $(x(0),y(0),z(0))\in [0,L]^3$.
It follows that for every point $(x(0),y(0),z(0))\in [0,L]^3$, we have
$(x(2),y(2),z(2))\in B_{(\sigma,\rho,K)}$.
\end{proof}

\begin{proof}Proof of Proposition~\ref{stablemanE}.
The time reversed system of \eqref{inequal} is a monotone dynamical system with respect to the partial
order induced by the octant $K=\{(x,y,z):x,y\ge 0,z\le 0\}$. So we make use of  Theorem 2.8 in \cite{S0}.
By irreducibility of the Jacobian matrix at $E$, there is an eigenvector $v\in K$ with nonzero components. Moreover, the unstable manifold of $E$
for the time-reversed system consists of $E$ and two strictly monotone trajectories $P_u\in E-K$ and $P_l\in E+K$ that are described in Theorem 2.8 in \cite{S0}. The monotone trajectory $P_u$ must exit $\mathbb{R}^3_+$ through the $y=0$ hyperplane because the derivative of $y$ is negative and bounded away from zero.
\end{proof}

\begin{proof}Proof of Theorem~\ref{persistence}.
We employ Theorem 8.17 in \cite{ST}, using the notation developed there. Let $p=(x,y,z)\in \mathbb{R}^3_+$ and $\rho(p)=\min\{x,y,z\}$.
$X_0=\{p \in \mathbb{R}^3_+: \rho(\phi(t,p))=0,\ t\ge 0\}$ consists of the stable manifold $S$ of $0$. The requisite compactness
assumption (H) of Theorem 8.17 follows from Theorem~\ref{attractor}. $\Omega=\cup_{p\in X_0}\omega(p)=\{0\}$ where $\omega(p)$
denotes the omega limit set of $p$. $M\equiv \{0\}$ is obviously an acyclic covering of $\Omega$ and it is an isolated invariant set in $\mathbb{R}^3_+$ because $0$ is a hyperbolic equilibrium (Hartman-Grobman Theorem). Finally, $M$ is weakly $\rho$-repelling, meaning that there is no
$p\in \mathbb{R}^3_+$ such that $\rho(p)>0$ and $\phi(t,p)\to M$ as $t\to \infty$. This follows because the stable manifold of $0$ is exactly
the $x=0$ facet $S$ where $\rho=0$. Consequently, we may apply Theorem 8.17 in \cite{ST} to conclude that $\phi$ is uniformly weakly $\rho$-persistent:
$\exists \eta>0$ such that $\rho(p)>0 \Rightarrow \limsup_{t\to\infty}\rho(\phi(t,p))>\eta$. We conclude uniform strong persistence, i.e., replacing limsup by liminf in the above implication, by invoking Theorem 4.5 in \cite{ST}.

The partition of the attractor $A=\{0\}\cup C\cup A_1$ is a consequence of Theorem 5.7 in \cite{ST}.
\end{proof}

\begin{proof}Proof of Theorem~\ref{solnbehavior}.
Part (a) was proved in \cite{WH2}.  Part (c) follows from
Theorems 1.1 and 1.2 of \cite{ZS} where we check that the determinant of the Jacobian at $E$ is $-kk_3k_5<0$.

The final assertion uses the uniform persistence established in Theorem~\ref{persistence}, which implies that
the set $D=[\eta,M]^3$ is an absorbing set for positive solutions of \eqref{scaled}.
Since $\frac{1}{T}\int_0^T \frac{x'}{x}ds=\frac{1}{T}\ln(x(T)/x(0))\to 0$
as $T\to\infty$ because $x$ is bounded and bounded away from zero, we see from the differential equation that $\lim_{T\to\infty}\frac{1}{T}\int_0^T (k-y(s))ds=0$. Indeed, this limit is uniform for all solutions starting in $D$ because $|\ln(x(T)/x(0))|$ is uniformly bounded on $D$.
Thus, $$\lim_{T\to\infty}\frac{1}{T}\int_0^T y(s)ds=k$$ uniformly for solutions in $D$ . The other limits follow more easily and they too are uniform in $D$.

Part (b) is more technical. We use the geometric approach of ruling out periodic orbits developed in \cite{LM} which uses the second additive compound  $Df(p)^{[2]}$ of the Jacobian matrix $Df(p)$, where $f(p)$ denotes the vector field \eqref{scaled}:
$$
Df(p)^{[2]}=\left(
              \begin{array}{ccc}
                k-y-k_3 & k_3 & 0 \\
                0 & k-y-k_5 & -x \\
                -k_5 & 0 & -(k_3+k_5) \\
              \end{array}
            \right)
$$
See, e.g., formula for $Df(p)^{[2]}$ in appendix of \cite{LSW}. As $k<k_3+k_5$, we may find $\epsilon\in (0,1/2)$ such that $k/(1-2\epsilon)<k_3+k_5$.
Consider the $3\times 3$ matrix-valued function $p\to A(p)$ where $A(p)$ is a diagonal matrix with diagonal entries $(1-2\epsilon)/x,\ (1-\epsilon)/x,\ -1/k_5$. Using the notation of \cite{LM}, we compute the matrix function $B(p)=A_f A^{-1}+ADf(p)^{[2]}A^{-1}$ where $A_f$ is the directional derivative of $A$ with respect to vector field $f$. We find that
\begin{equation}\label{B}
B(p)=\left(
       \begin{array}{ccc}
         -k_3 & k_3\frac{1-2\epsilon}{1-\epsilon} & 0 \\
       0 & -k_5 & k_5(1-\epsilon) \\
         \frac{x}{1-2\epsilon} & 0 & -(k_3+k_5) \\
       \end{array}
     \right)
\end{equation}
The Lozinski\u{i} measure $\mu(B)$ relative to the norm $|p|=\max\{|x|,|y|,|z|\}$ is (see pg 41 of \cite{C}):
$$
\mu(B)=\max\{-\epsilon k_3/(1-\epsilon), -\epsilon k_5, \frac{x}{1-2\epsilon}-(k_3+k_5)\}
$$
Now we use that  $\lim_{T\to\infty}\frac{1}{T}\int_0^T x(s)ds=k$ is uniform for all orbits starting in the absorbing set $B$ and that $k/(1-2\epsilon)<k_3+k_5$ to conclude
that $$\bar q_2=\limsup_{t\to\infty}\sup_{p\in D}\frac{1}{T}\int_0^T \mu(B(\phi(s,p)))ds<0.$$
Theorem 3.1 of \cite{LM} implies that there can be no nontrivial periodic orbit of \eqref{scaled}. By the Poincar\'e-Bendixson
Theorem for competitive three dimensional systems \cite{H0,HS,S}, $E$ attracts all solutions not starting on $S$.
\end{proof}


\end{document}